\documentclass[letterpaper,12pt]{article}
\usepackage[english,activeacute]{babel}
\usepackage[inner=2.5cm,outer=2.5cm]{geometry}
\usepackage{amssymb,amsmath,alltt}
\usepackage[utf8]{inputenc}
\usepackage{amsfonts}
\usepackage{graphicx}
\usepackage{multirow}
\newcommand{\bq}{\begin{equation}}
\newcommand{\eq}{\end{equation}}

\newcommand{\be}{\begin{eqnarray}}
\newcommand{\ee}{\end{eqnarray}}

\title{Self-consistent chaotic transport in a high-dimensional mean-field Hamiltonian map model}
\author{D. Martínez-del-Río$^\diamond$\footnote{dmr@mym.iimas.unam.mx},
       D. del-Castillo-Negrete$^\star$\footnote{delcastillod@ornl.gov},
       A. Olvera$^\diamond$\footnote{aoc@mym.iimas.unam.mx}, 
      R. Calleja$^\diamond\footnote{calleja@mym.iimas.unam.mx}$ 
      \\  \footnotesize{$^\diamond$ IIMAS-UNAM, Mexico, D.F. 04510} 
         \phantom{OOOOOOOOOOOOOOiiOOO}
        \\ \footnotesize{$^\star$ Oak Ridge National Laboratory, Oak Ridge, Tennessee 37831-8071}}
\date{}

\begin{document}

\maketitle

\begin{abstract}

Self-consistent chaotic transport is studied in a Hamiltonian mean-field model. The model 
provides a simplified description of transport in marginally stable systems including
vorticity mixing in strong shear flows and electron dynamics in plasmas.  Self-consistency is incorporated through a mean-field that couples all the degrees-of-freedom.  The model is formulated  
as a large set of $N$ coupled standard-like area-preserving twist maps in which the amplitude  and phase of the perturbation, rather than being constant like in the standard map, are dynamical variables. 
Of particular interest is the study of the impact of periodic orbits on the chaotic transport and coherent structures.
Numerical simulations show that self-consistency leads to the formation of a coherent macro-particle trapped around the 
elliptic fixed point of the system that appears together with an asymptotic periodic behavior of the mean field.  
To model this asymptotic state, we introduced a non-autonomous map  that allows a detailed study of the onset of 
global transport. A  turnstile-type transport mechanism that allows transport across instantaneous KAM invariant 
circles in non-autonomous systems is discussed. 
As a first step to understand transport, we study a special type of orbits referred to as sequential periodic orbits. Using symmetry properties we show that, through replication, high-dimensional sequential periodic orbits can be generated starting from low-dimensional periodic orbits.  We show that sequential periodic orbits in the self-consistent map can be continued from trivial (uncoupled) periodic orbits of standard-like maps using numerical and asymptotic methods. Normal forms are used to describe these orbits  and to find the values of the map parameters that guarantee their existence. 
Numerical simulations are used to verify the prediction from the asymptotic methods.
\end{abstract}

%|||||||||||||||||||||||||||||||||||||||||||||||||||||
\section{\large{Introduction}}
%|||||||||||||||||||||||||||||||||||||||||||||||||||||
There are many  problems in physics and the applied sciences where transport plays an
important role. Some examples include the 
dispersion of pollutants in the atmosphere and the oceans,  the efficiency of mixing in
chemical reactions, the magnetic confinement of fusion plasmas,  and  vorticity  mixing 
in fluids among others.  The study of transport
can be classified in  two groups: passive and  active
transport. In the case of passive transport, the elements which are
transported do not modify the velocity field of the flow while in the
active case, the advected field determines the velocity field through a self-consistent feedback mechanism.  

In the present paper, we focus on active transport.
The study is based on a simplified model of self-consistent transport in 
marginally stable systems including vorticity mixing in strong shear flows and electron dynamics in plasmas.
The model, originally presented in Ref.~\cite{del_castillo_2000} and later on studied in 
Ref.~\cite{Carbajal12},  consists of a large set of $N$ mean-field coupled standard-like area-preserving twist maps in which the amplitude  and phase of the perturbation (rather than being constant) are dynamical variables. An important part of our methodology is based on the use of normal form methods to describe the evolution of 
the dynamical variables of the system.

The rest of this paper is organized as follows.  Section 2 presents a brief summary 
of the model and some results originally contained in Ref.~\cite{del_castillo_2000} that will be used 
in the present work.
In Sec.~3, we study the asymptotic behavior of numerical simulations of 
the  self-consistent map model and in particular the formation of
coherent structures. To study these structures, we introduce a reduced transport 
map model that mimics  the asymptotic behavior and that can be interpreted as an uncoupled system in an external oscillatory field. Large iterations of this map are computed for particular sets 
of initial conditions to show that global diffusion can occur as a result of the 
oscillation of the perturbation parameter.
In Sec.~4, we study the properties of periodic orbits in the self--consistent 
map and introduce an asymptotic procedure to find a particular 
kind of periodic orbits that we name sequential periodic orbits.
In Sec.~5 we present a simple example illustrating how these 
sequential  periodic orbits can be computed using normal forms and present the 
relationships between the parameters of the map that should be satisfied 
for these orbits to exits. The values of 
the parameters of the  sequential periodic orbits obtained by numerical 
continuation methods are compared with the normal forms results.  It 
is also shown that these numerical and asymptotic results agree up to machine precision.
Section~6 present the conclusions.

%|||||||||||||||||||||||||||||||||||||||||||||||||||||
\section{\large{Hamiltonian mean field model of self-consistent transport}}
%|||||||||||||||||||||||||||||||||||||||||||||||||||||
The advection-diffusion equation
\begin{equation}
         \partial_t \zeta + {\bf{V}} \cdot \nabla {\zeta}= D\nabla^2 \zeta \, ,
         \label{addif_1}
  \end{equation}
where $\zeta$ is an advected scalar, $\bf{V}$ is the velocity field, and $D$ is the diffusivity is one of the fundamental models to study transport in the case of incompressible ($\nabla \cdot {\bf V}=0$) velocity fields. 
Restricting attention to $2$-dimensions, and introducing the stream-function $\psi$,
\begin{equation}
\label{def}
\mathbf{V}=\hat{\mathbf{z}}\times\nabla\psi \, ,
  \end{equation}
  Eq.~(\ref{addif_1}) can be written as
 \begin{equation}
         \partial_t \zeta - (\partial_y \psi)\,  \partial_x {\zeta}
         + (\partial_x \psi)\, \partial_y {\zeta} = D\nabla^2 \zeta \, .
         \label{addif_2}
  \end{equation} 
Broadly speaking, there are two different classes of transport problems. In the case of 
\emph{passive transport}, $\zeta$,  evolves without affecting the velocity field, $\mathbf{V}$. A typical example is the transport of low concentration pollutants in the atmosphere and the oceans.
In the case of \emph{active transport}, $\zeta$ determines $\mathbf{V}$ through a dynamical self-consistency constraint of the form $\mathcal{F}(\zeta,\mathbf{v})=0$ involving an integral and/or differential operator. A paradigmatic example  is the $2$-dimensional Navier-Stokes equation for an incompressible fluid. In this case, $\zeta$ in Eqs.~(\ref{addif_1}) and (\ref{addif_2}) denotes the fluid vorticity  that is self-consistently coupled to the velocity field according to the constrain 
\begin{equation}
\zeta=\left(\nabla\times\mathbf{V}\right) \cdot \hat{\mathbf{z}} = \nabla^2 \psi\, .
  \label{vort2}
\end{equation}

Active transport problems, are intrinsically nonlinear and as a result, harder to study than passive transport problems. 
For example, the well-known challenges in understanding fluid turbulence reside in the nonlinearity in  
Eqs.~(\ref{addif_2})-(\ref{vort2}).  It is thus of significant interest to develop simplified models that capture 
the basic elements of the self-consistent coupling within a relatively simple mathematical setting. 
One of these simplified descriptions is the Single Wave Model (SWM) originally proposed in plasma physics 
\cite{Oneil71,del_castillo_98} and used to study self-consistent transport in marginally stable fluids in the 
presence of strong shear flows \cite{del_castillo_2000}.
In general, the stream-function $\psi$ has a complicated spatio-temporal dependence. However, in the
SWM, $\psi$ is assumed to have the relatively simple form
\be
   \psi=-\frac{1}{2}y^2 + a(t)\,e^{ix} +a^*(t)e^{-ix}, \label{sw1}
 \ee
and the  vorticity stream-function self-consistent coupling in Eq.~(\ref{vort2}) reduces to 
 \be
   \frac{da}{dt} - iUa=\frac{i}{2 \pi} \int e^{-ix} dx  \int  \, \zeta(x,y,t) \, dy \, , \label{sw2}
 \ee
 where $a=a(t)$ is in general complex, and $U$ is a constant real parameter. 
 That is, whereas in the Navier-Stokes equation the  stream-function has a general spatio-temporal dependence that is 
determined at each instant of time by solving the Poisson equation in (\ref{vort2}), in the SWM 
 the spatial dependence of the stream-function is given and self-consistency only enters when one determines the 
temporal dependence of $\psi$ though the amplitude $a(t)$ according to the ordinary differential equation in 
(\ref{sw2}).  According to Eq.~(\ref{def}), in the SWM the velocity field consists of a linear component in the 
$x$-direction and a traveling-wave component in the $y$-direction, i.e. 
${\bf V}=y {\bf i} + {\rm Re} [i a(t) e^{ix}] \, {\bf j}$, where ${\rm Re} $ denotes the real part. 
 
Despite its relative mathematical simplicity, the SWM is able to capture important dynamics of the full Navier-Stokes equation. In particular, as Fig.~\ref{cateye_cont} from Ref.~\cite{del_castillo_2000} shows, the SWM  exhibit the standard Kelvin-Helmholtz instability leading to the formation of the familiar cats' eyes vorticity structure found in unstable shear flows.  The result in Fig.~\ref{cateye_cont}  was obtained from the direct numerical integration of Eqs.~(\ref{addif_1}) and (\ref{sw2}) with $D=0.001$, $U=-1$, and initial condition:
\begin{equation}
    \zeta(x,y,t=0)=e^{-y^2/2}\,[1-0.2\,y\,\cos(x)].
    \label{defect_z}
\end{equation}
Although for the sake of simplicity we have stressed the fluid-dynamics interpretation of the SWM, it is important to reiterate that as discussed in Refs.~\cite{del_castillo_98,del_castillo_2000,del_castillo_2002} the SWM has its origins in plasma physics and it can be equally applied to the study of electron dynamics in a Vlasov-Poisson plasma in a neutralizing ion background. In this case,   $\zeta$ corresponds to the single-particle electron distribution function, $\psi$ is the electrostatic potential, the vorticity equation becomes the Vlasov equation, 
and the $(x,y)$ variables correspond to the $(x,u)$-phase space variables. 

To  reformulate the SWM as a finite (but arbitrarily large) degrees-of-freedom Hamiltonian dynamical system, we assume from now on $D=0$ and  introduce the point-vortex representation
\be
   \zeta(x,y,t)=2 \pi \sum_{j=1}^{N}\Gamma_j\,\delta[x-x_j(t)]\,\delta[y-y_j(t)],
 \ee
 where $(x_k(t),y_k(t))$ denote the Lagrangian trajectories of the $k=1,2,...N$ point vortices with 
  intensities $\Gamma_k$. In this case,  it can be shown that Eqs.~(\ref{addif_2}), (\ref{sw1}) and (\ref{sw2}) imply the following set of Hamiltonian differential equations determining $(x_k(t),y_k(t))$
 \begin{eqnarray}
    &&\frac{dx_k}{dt}= \frac{\partial H}{\partial y_k} \, , \qquad \frac{dy_k}{dt}= -\frac{\partial H}{\partial x_k} ,
     \label{xy_dis} 
  \end{eqnarray}
  where the Hamiltonian is
  \begin{equation}
H= \sum_{j=1}^{N} \left[\frac{1}{2}y_j^2 -  a(t) e^{i x_j} - a^*(t) e^{-i x_j}\right] \, ,
\end{equation}
and the function $a=a(t)$ is determined from 
 \begin{eqnarray}       
   &&\frac{da}{dt} - iUa= \frac{i}{N} \sum_{j=1}^{N} \Gamma_j\,e^{-ix_j} \label{a_dis} \, ,
  \end{eqnarray}
  which is the point-vortex representation of the SWM self-consistent vorticity stream-function relation in 
  Eq.~({\ref{sw2}).
  The Hamiltonian model in Eq.~(\ref{xy_dis}) is a mean field model in the sense that the dynamics of  the particles $(x_k,y_k)$ is determined by the time dependent field $a(t)$ which according to Eq.~(\ref{a_dis}) depends on the mean properties of the particles' positions.  

Defining
\begin{equation}
  a= \sqrt{J}\,e^{-i\theta}, \qquad p_k=\Gamma_k\,y_k,
\end{equation}
equations (\ref{xy_dis}) and (\ref{a_dis}) can be written as a full $N+1$-degrees of freedom Hamiltonian system
 \begin{eqnarray}
    \frac{dx_k}{dt}= \frac{\partial \mathcal{H}}{\partial p_k}, \;\;\;&&
    \frac{dp_k}{dt}= - \frac{\partial \mathcal{H}}{\partial x_k} ,
     \label{xy_ham} \\
   \frac{d \theta}{dt}= \frac{\partial \mathcal{H}}{\partial J}, \;\;\;&&
    \phantom{i}\frac{dJ}{dt}= - \frac{\partial \mathcal{H}}{\partial \theta} . \label{jt_ham}
  \end{eqnarray}
where,
\begin{equation}
    \mathcal{H}= \sum_{j=1}^{N} \left[\frac{1}{2\Gamma_j}p_j^2 -2\Gamma_j\sqrt{J}\cos (x_j - \theta) \right] - UJ.
\end{equation}
%////////////////////////////
 \begin{figure}[h!]
     \centering
     \includegraphics[height=6.8cm]{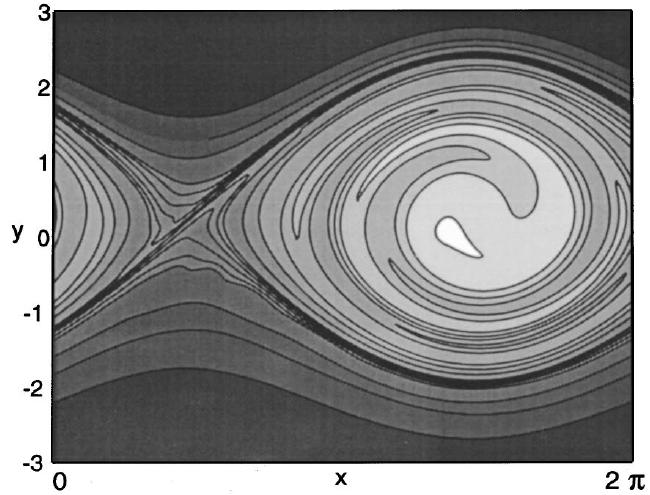} 
     \caption{\footnotesize{Cat's eye formation and vorticity mixing in the single wave model, Eqs.~(\ref{addif_2}),  (\ref{sw1}) and (\ref{sw2}), with initial conditions in Eq.~(\ref{defect_z}).  The gray scale denotes the vorticity values with white corresponding to $\zeta=1$ and dark gray corresponding to $\zeta=0$} 
     [After Ref.~\cite{del_castillo_2000}].}
     \label{cateye_cont}
 \end{figure}
%///////////////////////////

The last step to construct the map model is to perform
a  time discretization of Eqs.~ (\ref{xy_ham}) and  (\ref{jt_ham}) to obtain the self-consistent standard map model originally proposed in Ref.~\cite{del_castillo_2000},
\begin{subequations}
\label{mapaut}
\begin{align}
  x_{k}^{n+1}&=x_{k}^{n} + y_{k}^{n+1} \label{mapaut1}\\
  y_{k}^{n+1}&=y_{k}^{n} - \kappa^{n+1}\, \sin(x_{k}^{n}-\theta^{n}) \label{mapaut2}\\
  \kappa^{n+1}&= \sqrt{( \kappa^{n})^{2}+ (\eta^{n})^2} + \eta^{n} \label{mapaut3}\\
  \theta^{n+1}&= \theta^{n} - \Omega + {\frac{1}{\kappa^{n+1}} \frac{\partial \eta^n}{\partial \theta^{n}}}
 \label{mapaut4}
\end{align}
\end{subequations}
where $k= 1,2,\dots,N$, $x_k^n,\theta_n \in [0,2\,\pi)$, $ y_k^n, \kappa^n \in {\mathbb R}$, $\Omega= U \tau$ and
the $\eta^n$ is defined as:
\begin{equation}
\eta^n:=\sum_{k=1}^{N}\,\gamma_{k}\,\sin(x_{k}^{n}-\theta^{n}),
 \label{mapaut5}
\end{equation}
where $\tau$ is the discretization time step.
The map in Eqs.~(\ref{mapaut}) is a $2N+2$ dimensional map with $N+1$ parameters. The $N$ parameters,
$\gamma_k$, represent the intensities of the point vortices, and the constant parameter $\Omega$ is related to the parameter $U$  in the SWM. 

Note that Eqs.~(\ref{mapaut1})-(\ref{mapaut2}) have the structure of the well-known symplectic standard map 
\begin{subequations}
\label{stnmap}
\begin{align}
  x^{n+1}&=x^{n} + y^{n+1} \label{stn1}\\
  y^{n+1}&=y^{n} - \varepsilon\, \sin(x^{n} -\phi) \label{stn2}
\end{align}
\end{subequations}
with constant, fixed perturbation $\varepsilon$, and phase $\phi$ (which is usually taken equal to zero). 
However, in the self-consistent map, the perturbation parameter, 
$\kappa^n$, as well as the phase, $\theta^n$, depend on the iteration, $n$, and their dynamics is dictated by 
 Eqs.~(\ref{mapaut3})-(\ref{mapaut4}) which are the discrete map version of the self-consistent SWM coupling 
 in Eq.~(\ref{a_dis}). Although  Eqs.~(\ref{mapaut3})-(\ref{mapaut4}) do not define a symplectic transformation, they can be rewritten as a symplectic (although implicit) map in the mean field degrees of 
 freedom \cite{del_castillo_2000}.
The analogy with the standard map allows us to interpret the map in Eqs.~(\ref{mapaut1})-(\ref{mapaut2}) as $N$   coupled  oscillators through their phase and amplitudes by the  \emph{mean field} in Eqs.~(\ref{mapaut3})-(\ref{mapaut4}). We will take advantage of these similarities not only to identify the equations of the map, but also in its perturbation analysis.
Note that in the limit $\gamma_k\rightarrow0$ in Eq.~(\ref{mapaut5}), the   oscillators in (\ref{mapaut})  depending on $\theta^0$ and $\Omega$ decoupled and their equations simply correspond to $N$ copies of the standard map. 

%|||||||||||||||||||||||||||||||||||||||||||||||||||||
\section{\large{Coherent structures and transport}}
%|||||||||||||||||||||||||||||||||||||||||||||||||||||
The evolution of the self-consistent map in Eqs.~(\ref{mapaut}) has been studied
for different initial conditions. Figures~\ref{catseye_map} and \ref{evol} show the results of a simulation 
of $N=13,440$ coupled maps with
initial conditions $\{ (x_k^0, y_k^0) \}$ uniformly distributed on the region  $[0,2\pi]\times[-0.3,0.3]$ in the $(x,y)$ plane and 
$\gamma_k=3\times 10^{-6}$ for $k=1, \ldots N$. The initial condition of the mean-field was
 $\kappa^0=10^{-4}$ and $\theta^0=0$,  and we assumed $\Omega=0$.

 The simulations  show that the self-consistent map reproduces the coherent structures observed in the single wave model (Figure \ref{cateye_cont}).
In particular, while a group of particles exhibit coherent behavior trapped in the center of the cat's eye, 
those particles located in the separatrix exhibit a strong dispersion.

%////////////////////////////
  \begin{figure}[h!]
   \centering
     \includegraphics[height=3.2cm]{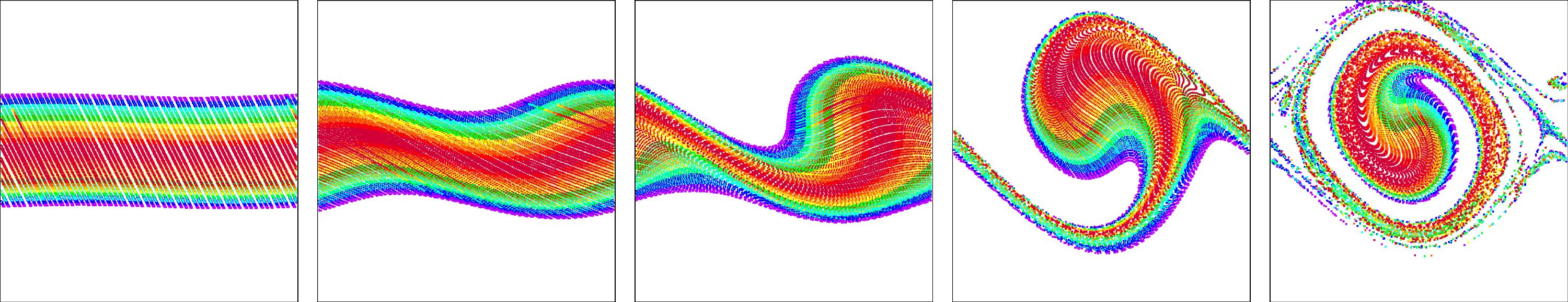}
     \caption{\footnotesize{Time evolution of the self-consistent map in Eq.~(\ref{mapaut}) with $N=13440$ initial conditions uniformly distributed in $[0,2\pi]\times[-0.3,0.3]$ with $\kappa^0=10^{-4}$, $\gamma_k=3\times 10^{-6}$, $\theta^0=0$ and $\Omega=0$. 
     The frames show the instantaneous coordinates of the $N$ initial conditions, at $n=2, 6, 12, 20$ and $66$
     in the region $[0,2\pi]\times[-0.8,0.8]$ of the $(x,y)$ plane.
     The colors label the $y$-coordinate of the initial condition,  with red denoting $y_k^0$ close to $y=0$ and blue further away.}}
      \label{catseye_map}
  \end{figure}
%////////////////////////////
  
The evolution of the mean field is represented by the variables
$(\kappa^n,\theta^n)$, which are shown in Fig.~\ref{evol}. The behavior of
$\kappa^n$ starts with a fast growth  until it achieves a maximum value, after that,
the $\kappa^n$ oscillates around a mean value $\bar{\kappa}$ and the amplitude of
oscillation is bounded by $\Delta\,\kappa$.  A similar situation is observed with
the behavior of $\vartheta^{n+1}= \theta^{n+1}-\theta^{n}$. 
Different types of dynamics, including cases in which the mean field decays to zero or saturates at a constant fixed value can be found in Refs.~\cite{Bofetta03, Carbajal12} for similar self-consistent maps. 
Note that in Fig.~\ref{evol}, $\kappa^n$
does not reach the critical value $\kappa_c=0.971635406$,
which is  the value of the
parameter $\kappa$ which corresponds to the destruction of the invariant
circle with rotation number equal to the \emph{golden mean}\footnote{The last invariant circle not homotopic to a point. It must be noted that due the choice of scale of the map, the rotation number is $\omega=2\pi\gamma$, instead of the golden mean: $\gamma=\frac{\sqrt{5}-1}{2}$.} of the standard map \cite{Greene, Meiss}. For  any $\kappa < \kappa_c$ there is no
 global diffusion in the standard map because of the existence of invariant
circles,  that give rise to  barriers in phase space. 
The existence or non-existence of global diffusion in the self-consistent map depends in a 
nontrivial way on the dynamics of $\kappa^n$.
On a more fundamental level, the observed rapid growth of $\kappa^n$ for a given initial condition is closely related to 
the linear stability properties of the corresponding initial condition in the single wave model.  In particular, 
Ref.~\cite{del_castillo_98} presents the necessary and sufficient conditions for the linear stability (i.e., exponential 
growth of the mean field amplitude) of a given initial condition in the context of the continuous, $N\rightarrow \infty$
limit. These ideas might help understand the conditions for the growth of  
$\kappa^n$. However, one must be careful before drawing conclusions as the self-consistent map model discussed here 
is obtained by simplifying drastically the single wave model by approximating the continuous limit when  
$N\rightarrow \infty$.

%////////////////////////////
 \begin{figure}[h!]
 \centering
  \includegraphics[height=6.5cm]{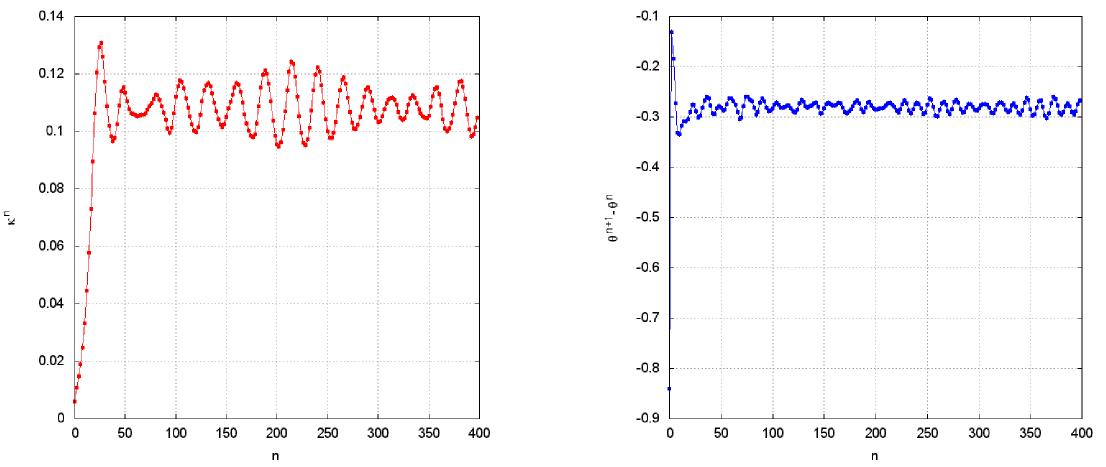}
  \hspace{-0.4cm}
  \begin{minipage}[b]{0.55\linewidth}
     \centering
     \phantom{OOOOOO}(a)
  \end{minipage}
  \begin{minipage}[b]{0.42\linewidth}
     \centering
     \phantom{O}(b)
  \end{minipage}
    \caption{\footnotesize{{Time evolution of the mean-field variables in the self-consistent map in Eq.~(\ref{mapaut}) with initial conditions taken in a uniform grid in $[0,2\pi]\times[-0.3,0.3]$ with $N=13440$, $\kappa^0=10^{-4}$, $\gamma_k=3\times 10^{-6}$, $\theta^0=0$ and $\Omega=0$.   
    $\kappa^n$ is shown in (a) and $\vartheta^{n+1}=\theta^{n+1} - \theta^n$ in (b). }}}
  \label{evol}
 \end{figure}
%////////////////////////////
%////////////////////////////
\begin{figure}[h!]
     \centering
     \includegraphics[height=9cm]{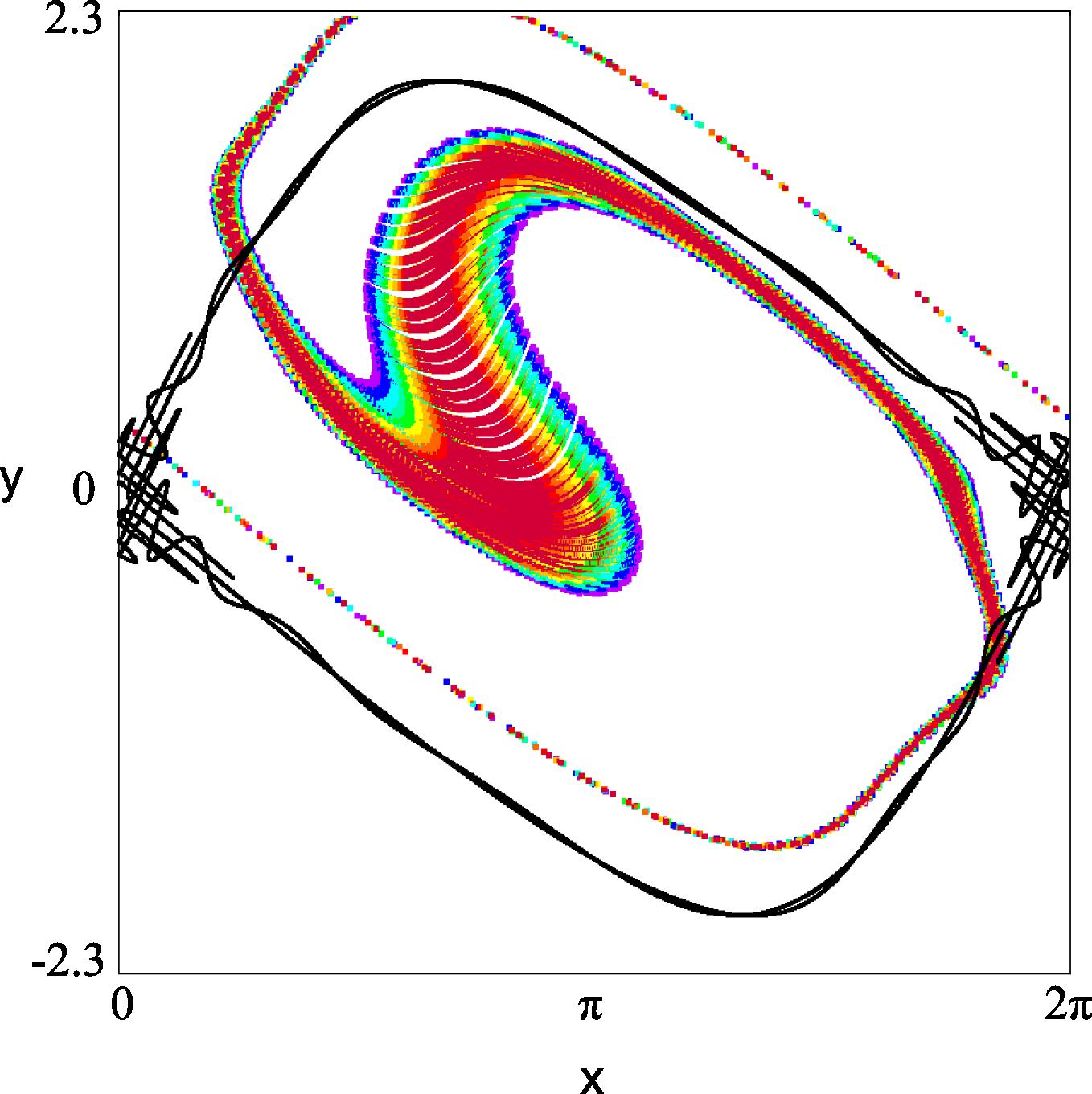}    
  \caption{\footnotesize{Plot of the projected phase space coordinates $(x_i,y_i)$ on the $(x,y)$-plane 
 of a simulation of the self-consistent map in Eq.~(\ref{mapaut}) in the oscillatory $\kappa$ regimen.
 Also shown in black is the heteroclinic tangle generated by the unstable invariant 
manifold of the hyperbolic fixed point of the standard map with perturbation parameter equal to $\kappa^n$. 
The initial conditions are the same as those in Fig.~\ref{catseye_map}, except that to enhance the heteroclinic tangle the higher values $\kappa^0=0.005$ and $\gamma=0.0005$ were used. }}
\label{ojo}
\end{figure}
%//////////////////////////// 

Perhaps contrary to the intuition, it is observed that global diffusion exists 
even when $\bar{\kappa} < \kappa_c$. 
It is also worth mentioning that when the \emph{instantaneous} coordinates $(x_k^n,y_k^n)$ of each
degree-of-freedom are plotted on the same plane like in Fig.~\ref{catseye_map}, the amplitude and shape of the  cat's eye structure is in 
good agreement with the invariant manifolds emanating from the hyperbolic fixed point of the standard 
map calculated with a perturbation parameter equal to the \emph{instantaneous} value of 
$\kappa^{n+1}$.
This gives rise to the following question:  What is the mechanism
that allows the diffusion across the invariant curves on the self-consistent map? 
In Figure~\ref{ojo} we observe the formation, for relatively small times, of the macro particle coherent structure 
trapped around the elliptic 
fixed point, and at the same time we have the formation of the heteroclinic tangle 
responsible for the high mixing region around 
the separatrix of the cat's eye. 
 
%////////////////////////////
  \begin{figure}[h!]
     \centering
     \includegraphics[height=10.1cm]{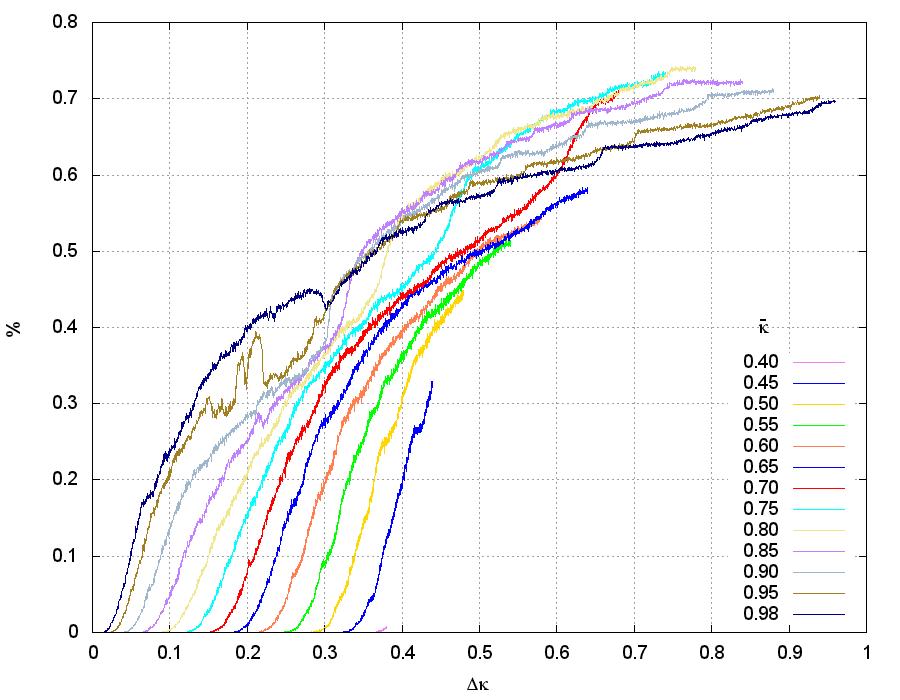} 
     \caption{\footnotesize{Onset of global transport in the non-autonomous map (\ref{map}) for different values of ${\bar{\kappa}=(\kappa_{1}+\kappa_{2})/2}$. 
     The $x$-axis corresponds to the difference ${\Delta\kappa=\kappa_{2}-\kappa_{1}}$ { for a fixed value of $\bar{\kappa}$.} 
     The $y$-axis corresponds to the percentage of  initial conditions, $(x_k^0,y_k^0)$, taken from a $100\times100$ regular 
     partition of the rectangle $[0,2 \pi]\times[0, {\pi/5}]$, that escape due to global transport in the sense that  after $M=100,000$ iterations satisfy 
     $y^M_k\geq \pi/5 + 2\pi$.}}
     \label{zfig}
 \end{figure}
%////////////////////////////

To address this question, and motivated by the observed oscillation of $\kappa^n$, we propose the 
following non-autonomous map consisting of two copies of the 
standard map
applied sequentially with alternating values of $\kappa$,

\begin{equation}
  \begin{array}{rcl}
  y^{n+1} & = & y^n - \kappa^n \sin x^n \\
          &  &  \\
  x^{n+1} & = & x^n + y^{n+1} \\
  \end{array} 
   \phantom{OOO} \mathrm{where} \phantom{OO} 
  \kappa^n = \left\lbrace \begin{array}{cl}
\kappa_1 & \;\;\;\; {\rm if} \; n \; {\rm is}\;{\rm odd} \\
  & \\
\kappa_2 & \;\;\;\; {\rm if} \; n \; {\rm is}\;{\rm even} \\
\end{array} \right.
\label{map}
\end{equation} 
Without loss of generality, we will assume that $\kappa_2$ is greater than $\kappa_1$.

To explore the onset of global diffusion for values of $\kappa_1$ and $\kappa_2$ less than $\kappa_c$, we considered a  set of $N$  initial 
conditions uniformly distributed in the region 
$[0,2\pi]\times[0,y_{\rm max}]$  
{ and the map is iterated $n$ times, with $n$ less than some given maximal value $M$}.
We then found the number of initial conditions that reached the semi-space $y > y_{\rm max} + 2\pi $. 
In turn, this means that one orbit could pass through the invariant circles which exist in
the standard map with  $\kappa_1$ and $\kappa_2$ less than $\kappa_c$. Figure \ref{zfig} shows 
the percentage of initial orbits that reached the semi-space $y > y_{\rm max} + 2\pi$
for a given $\kappa_1$ and $\kappa_2$. 
For convenience we show the results as function of $\Delta\kappa=\kappa_2-\kappa_1$ and $\bar{\kappa}=\frac{\kappa_1+\kappa_2}{2}$. Note that in most cases $\bar{\kappa}$ is lower 
than $\kappa_{c}$ and yet there are cases of global diffusion even when $\kappa_1<\kappa_2<\kappa_c$. 
This results provide numerical evidence of the existence of global diffusion. However, as shown in 
Fig.~\ref{zfig} there are also cases with no diffusion, a result consistent with Moser's theorem 
\cite{Moser62} that ensures that for values of $\kappa_1$ and $\kappa_2$ sufficiently small and for which 
the twist condition is satisfied, the map (\ref{map})  should have invariant circles that forbid global transport. 

%////////////////////////////
\begin{figure}[h!]
 \centering
 \includegraphics[height=8.5cm]{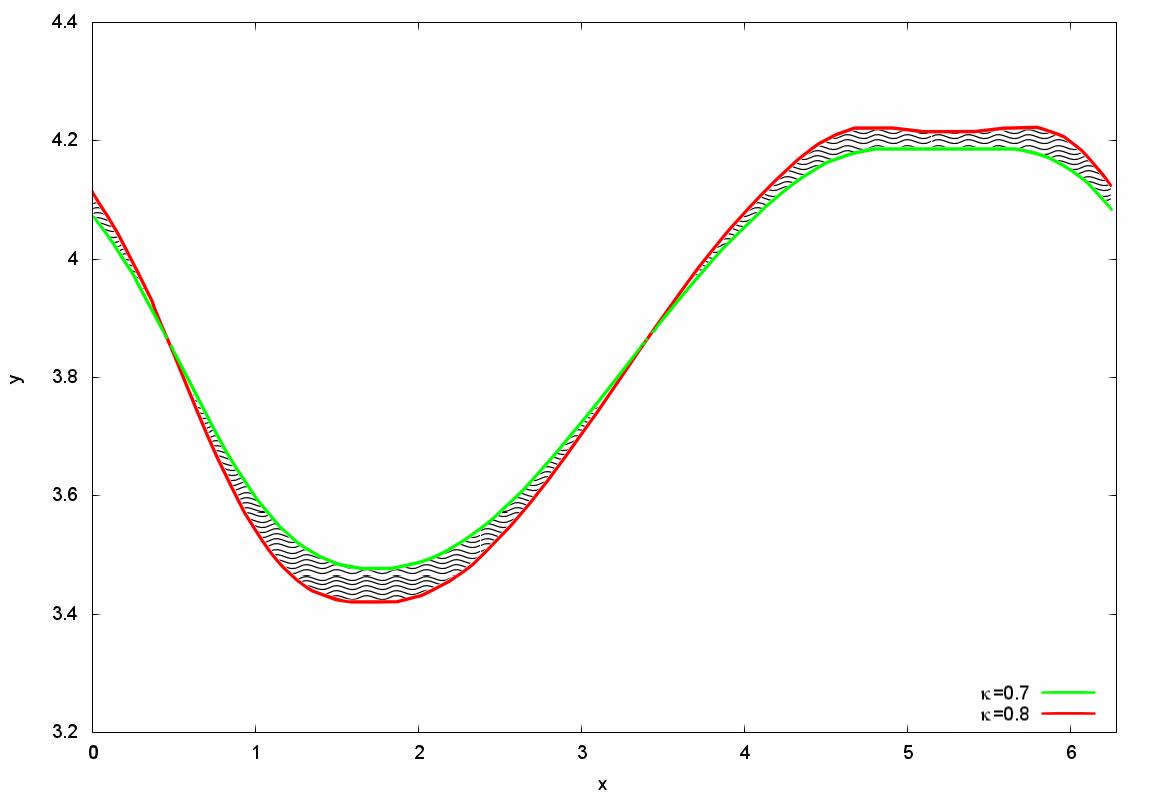}    
 \caption{\footnotesize{Illustration of the turnstile-type transport mechanism across invariant circles in the non-autonomous standard map. The plot shows 
the invariant circles (approximated using the parametrization method
in Ref.~\cite{Calleja09}.) with golden mean rotation number in the 
standard map (\ref{stnmap}) for two different values of the perturbation parameter 
$\kappa_1=0.7$  and $\kappa_2=0.8$.
Due to the periodic switching between these two values of $\kappa^n$  in the non-autonomous standard map in Eq.~(\ref{map}), points in the shaded region can cross the golden mean tori even though the perturbation is below the critical value 
$\kappa_c=0.971635406$  
for the onset of global transport in the standard map.
}}
\label{turnstyle_fig}
\end{figure}
%//////////////////////////// 

According to Birkhoff's theorem, an invariant
circle is the graph of a function $y=f(x)$ and  the position of this curve in the
$(x,y)$ plane is determined by the value of $\kappa$ \cite{Meiss}. As shown in Fig.~\ref{turnstyle_fig}, the
shape and position of invariant curves with a given rotation number for  $\kappa_1\neq 
\kappa_2$ in general do not coincide, and the set of points located between these two curves
can move upward and downward when we iterate the map (\ref{map}). Like in a revolving door,  every two iterations, points between the two invariant circles move across the invariant circles.
This turnstile-like type mechanism does not appear for all $\kappa_1$ and
$\kappa_2$.  In fact, we can appreciate in Fig.~\ref{zfig} that the probability to jump
across the invariant circles goes to zero when $\Delta \, \kappa$ approaches zero. 
This means that there exist invariant circles for the composed map
(\ref{map}) for small values of $\Delta\,\kappa$. 

Note that even if two standard maps with perturbation amplitudes $\kappa_1$ and $\kappa_2$ do not have 
common invariant curves, there is no reason to assume that their composition has no invariant curves.
To determine the threshold of $(\kappa_1, \kappa_2)$ values for which  
the map (\ref{map}) exhibits global diffusion more simulations were performed for a
wider range of values. 
{ To do this, we used the following criterion: if there is an initial condition 
{ $y_0 \in [0,y_m]$, with 
$y_m \in (0,2\pi)$, for which  
$y^n > y^*=y_m+2\pi \ell$ or 
$y^n \leq y^{**}=-2\pi \ell$, for some $\ell \in \mathbb{N}$}
for a given $n$, then there is global transport, i.e. the system has no invariant circles. Figure~\ref{cuerno} shows the minimum value of $\kappa_2$ for a given $\kappa_1$ such that the simulations satisfy this criterion,  for $\ell = 1, 2$. Note that for the case $\ell=2$ we have increased the number of iterations to obtain estimates comparable to $\ell=1$. We emphasize that this criterion assumes that the map is invariant under the transformation $(x_n,y_n) \rightarrow (x_n, y_n + 2 \pi)$ 
 which is the case of the map under study}.
We found a symmetrical bifurcation diagram in positive quadrant 
of the $(\kappa_1, \kappa_2)$ parameter space. Points inside the  ``horn" correspond to 
cases were none of the  initial conditions reached $y*$. In other words, there is not 
global diffusion for values of $\kappa_1$ and $\kappa_2$ inside the  horn.

{ It can also be observed in Fig.~\ref{cuerno}  that the numerically determined thresholds for global transport for parameter values $(\kappa_1,0)$ and $(0,\kappa_2)$ satisfy  to a very good approximation $\kappa^*_1=\kappa^*_2= \kappa_c/2$. This follows directly from the fact that, for the case $(\kappa_1,0)$, the non-autonomous standard map for every other iteration reduces to
}
    \begin{subequations}
    \label{mapk1v1_old}
    \begin{align}
            x_{n+1}&= x_n + 2y_n - 2\kappa_1               
                       \sin( x_n) \label{mapk11v1}\\
            y_{n+1}&= y_n - \kappa_1 \sin( x_n)   \label{mapk12v1} 
    \end{align}
    \end{subequations}  
{ Upon the change of coordinates, $\{X=x,Y=2y\}$, this map becomes the standard map (\ref{stnmap}) with $\varepsilon=2\kappa_1$, which as it is well known exhibits global chaos for $\varepsilon=\kappa_c$.  An analogous reduction happens for the case $(0,\kappa_2)$.}          
Further details of these calculations will be reported in a forthcoming publication
where we also study negative 
$(\kappa_1, \kappa_2)$  values and  
 examine the rotation numbers of the most robust invariant circles. 

%////////////////////////////
\begin{figure}[h!]
\centering
 \includegraphics[height=8.5cm]{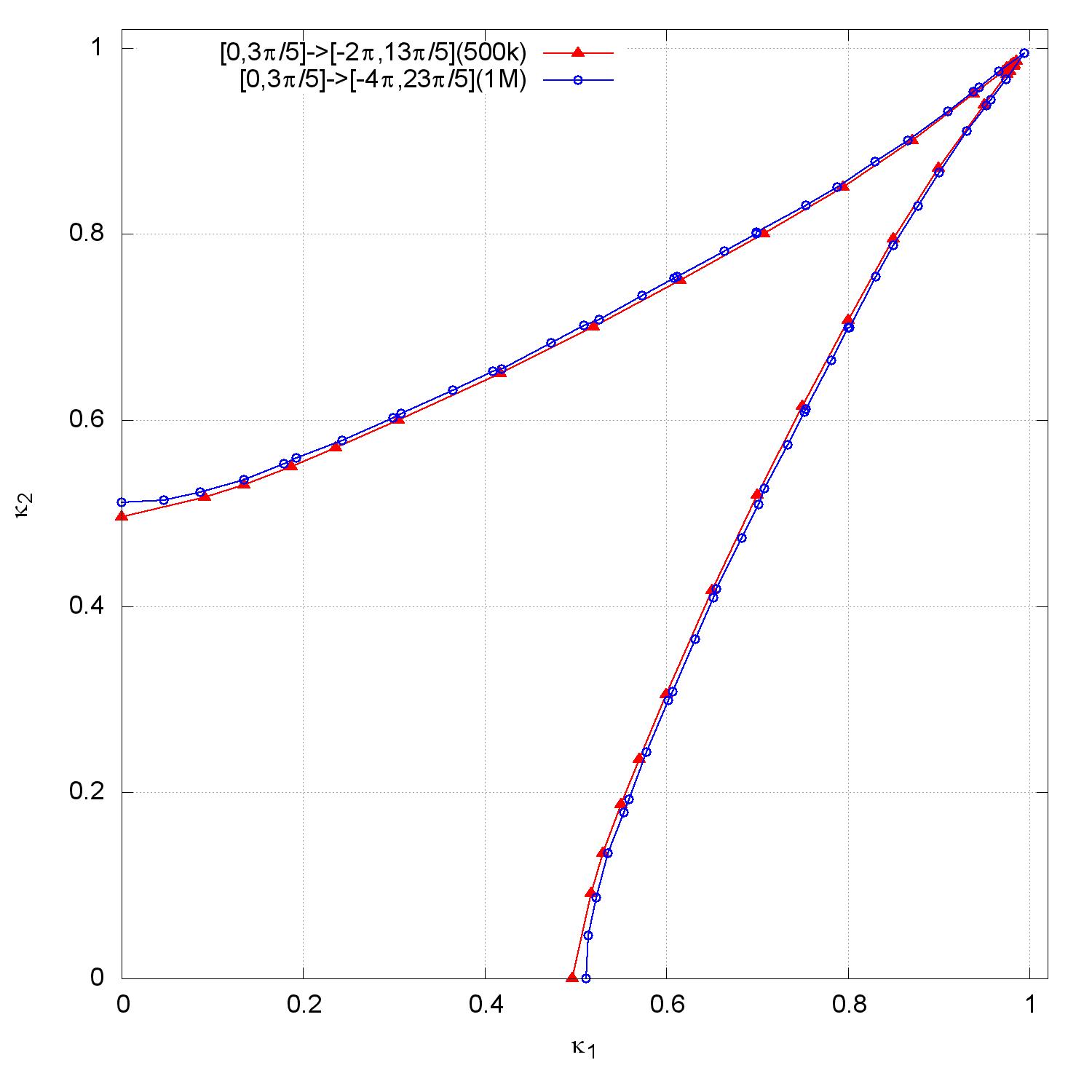}  
 \caption{\footnotesize{Critical values $(\kappa_1, \kappa_2)$ for global transport in the non-autonomous  standard map in Eq.~(\ref{map}). The area  outside the ``horn" corresponds to parameter values for which particles escape in the sense that at least one of $10^4$ initial condition uniformly distributed on the rectangle $[0,2 \pi ]\times[0, 3\pi/5]$ exhibited a displacement with: (i)
  ${y> 3\pi/5 + 2\pi}$ or ${y< - 2\pi}$  after $5\times 10^5$ iterations (red) and (ii) ${y>  3\pi/5 + 4\pi}$ {or} ${y<  -4\pi}$ 
     after $10^6$ iterations of the map (blue).}}
     \label{cuerno}
\end{figure}
%//////////////////////////// 

%|||||||||||||||||||||||||||||||||||||||||||||||||||||
\section{\large{Periodic orbits}}
%|||||||||||||||||||||||||||||||||||||||||||||||||||||
As a first step to understanding transport 
and the formation of coherent structures we present in this section a study of the periodic orbits of the self-consistent map. 
In particular, we determine the asymptotic value
of $\bar{\kappa}$ and $\Delta\,\kappa$ for an special kind of initial conditions 
and parameter values. To do this, the specific form of the mean field auxiliary 
variable $\eta^n$ (\ref{mapaut5}) is taken into account to try to determine which are the 
sets of initial conditions that correspond to $\eta^n$'s that are  as small as possible. 

Finding periodic orbits for high dimensional maps
is a complex problem and numerical simulation is in principle
the only procedure to estimate the asymptotic behavior of the mean field
variables. However, when using simulation it is hard 
to predict the result from the initial conditions and usually a large number 
of iterations are needed. 
 In our case, the dimension of the self-consistent map (\ref{mapaut})
has to be rather big  (of order $10^4$) so that the system can achieve similar 
patterns to that of
the single-wave model. In practice it is not possible to approach this problem
using analytic or asymptotic tools. One possibility is to reduce the complexity
of the problem by only taking into account  the periodic solutions of the map
(\ref{mapaut}). 
 In general, the structure of periodic solutions determines the
behavior of any dynamical system, our goal is to show that  periodic orbits
are closely related  to the mean-field variables { (in particular $\bar{\kappa}$ and $\Delta\,\kappa$) and the choice of values of the parameters.

In the case of low-dimensional maps, the use of symmetries can greatly
simplify the search of periodic orbits \cite{Kook}.  For the self-consistent
map, we can simplify the problem by assuming  $\gamma_k = \gamma$ constant,
for $k=1,\dots,N$ in (\ref{mapaut3}) which implies the following properties:

\begin{enumerate}

  \item Let ${\bf z}$ be a periodic orbit of map (\ref{mapaut}) 
  with dimension ${2N+2}$.  If we set the  initial condition of this orbit to be,

\begin{equation}
{\bf z}^0=(x_1^0,y_1^0,x_2^0,y_2^0,\dots,x_N^0,y_N^0,\kappa^0,\theta^0)
\;\;\;,
\label{orbN}
\end{equation}
then any permutation of the pairs $(x_i^0,y_i^0)$, for $i=1,\dots,N$, produces
a new periodic orbit with the same period. This is because the term that couples
the standard maps only depends on the average of the variables $x_i$.

  \item For any periodic orbit ${\bf z}$ of the  map (\ref{mapaut})
   with dimension
${2N+2}$, we can increase the dimension of the solution by replicating the orbit ${\bf
z}$ and reducing the strength of $\gamma$ by one  half. Thus, for
any ${\bf z}$ given in (\ref{orbN}), we can generate a new periodic solution of
dimension $2 N$ of the form,
$$
{\bf z}^0=(x_1^0,y_1^0,x_2^0,y_2^0,\dots,x_N^0,y_N^0,
          x_1^0,y_1^0,x_2^0,y_2^0,\dots,x_N^0,y_N^0,\kappa^0,\theta^0)
\;\;\;,
$$
with the strength of vorticity equal to $\gamma/2$. The dimension of this new  orbit is ${4N+2}$ and the period of the orbits is the same as (\ref{orbN}). In this case,
since $\gamma$ is rescaled by a factor of $1/2$, the function $\eta$ in (\ref{mapaut3})
preserves the value given when we calculate it at the orbit (\ref{orbN}) since we have to sum 
$N$ terms twice.

\end{enumerate}

  With these two properties, we can generate a large set of periodic
orbits: starting with low dimensional map corresponding to $N=1$ and a strength
parameter $\gamma$, we can find periodic orbits of period $M$. Using
the property 2, the orbit can be replicated $s$ times, such that we
get an orbit of dimension $2^s$ but with reduced strength in the parameter
$\gamma'=\gamma/2^s$.

 { Another possibility is to construct a periodic orbit by imposing that the 
iteration of each  initial condition almost coincides with the previous point in the orbit
of the next initial condition. This is, if $(x_i^{n},y_i^{n})$ is the $i$-th initial condition
and $(x_i^{n+1},y_i^{n+1})$ is the next point in the periodic
orbit, then $(x_i^{n+1},y_i^{n+1}) \approx (x_{i+1}^{n},y_{i+1}^{n})$. 
We will call this kind of periodic orbits as 
\emph{sequential periodic orbits} (SPO's).} To generate
this kind of orbits,  we compute a periodic orbit of period
$\tau$ of the standard map for a fixed value of $\kappa$ and
$\theta$. The iterations of this orbit give the following set of pairs
of coordinates: 
\bq
(x_0,y_0) \;\; (x_1,y_1) \;\; \dots \;\;
(x_{\tau-1},y_{\tau-1})\;\;\;\;. \label{orbit}
\eq
With this sequence of points,  we construct a periodic orbit of 
dimension ${2\tau+2}$ with  period $\tau$ for the  map (\ref{mapaut}),
with initial condition,
$$
(x_0,y_0,x_1,y_1,x_2,y_2, \dots , x_{\tau-1},y_{\tau-1},\kappa,\theta)
\;\;\;.
$$
For $\gamma=0$, this orbit follows a sequential pattern, that is,
the coordinates $(x_i,y_i)$ in the position $i$ are shifted to the
position $i+1$ for each iteration of the map.
It is clear that this is not a periodic orbit for $\gamma \neq 0$
because the mean-field variables, $\kappa$ and $\theta$,  are going to
change while the orbit is iterated. Nevertheless, this is a good guess
for finding a periodic orbit using numerical or asymptotical procedures
for  small values of $\gamma$.

Based on the previous ideas, our numerical procedure to find a periodic orbit of (\ref{mapaut}) with period $\tau$ and  dimension ${2\tau+2}$ consists of the following steps:

\begin{enumerate}
  %!!!!!!!!!!!!!!!!!!!!!!
  \item For given  values of $\theta$ and $\kappa$, we  find a symmetric
periodic orbit of period $\tau$ of the standard map using the symmetry
lines of the map \cite{Greene}.

  %!!!!!!!!!!!!!!!!!!!!!!
  \item With this orbit, a sequential periodic orbit (\ref{orbit}) is formed. For $\gamma=0$ this is a periodic orbit of the self-consistent map (\ref{mapaut}).
  
  %!!!!!!!!!!!!!!!!!!!!!!
  \item Using a continuation method, the sequential {periodic} orbit is
computed for small values of $\gamma$ \cite{Doedel} and continued for large values of $\gamma$. 

  %!!!!!!!!!!!!!!!!!!!!!!
  \item The dimension of the orbit is then increased using the
replication property of the map (\ref{mapaut}). %\ref{automap}.
 For each replication, the strength parameter $\gamma$ is divided by two.

  %!!!!!!!!!!!!!!!!!!!!!!
  \item { The value of the auxiliary variable $\eta^n$ is small in all iterations due the almost equidistribution on the angular variables $x_i^n$.}
\end{enumerate}

Following this procedure, we can compute periodic solutions for small periods. The convergence  
of the continuation--method requires adjusting the value of the parameter $\Omega$ in 
(\ref{mapaut4}), otherwise the variable $\theta$ might not converge to a value. 
After some numerical experiments, we concluded that for the method to converge,  $\Omega$ 
should be a function of $\kappa^0$. Therefore, we have included $\Omega$ in the numerical 
method as a free parameter which has to be determined by the continuation process. 

We proceed in a similar way to determine periodic orbits using a perturbative method. The idea is to compute a sequential periodic orbit in an asymptotic way. We have to begin the process with an integrable map, in our case this map is (\ref{mapaut}) when the parameter $\gamma$ and the variable $\kappa^n$ are zero. For this integrable case, the periodic solution with period $q$ is:
\begin{equation}
\left( x_1^0=2\pi\frac{p}{q},y_1^0=2\pi\frac{p}{q},
       x_2^0=2\pi\frac{2p}{q},y_2^0=2\pi\frac{p}{q},...,
       x_q^0=2\pi\frac{qp}{q},y_q^0=2\pi\frac{p}{q},
       \kappa^0=0,\theta^0=\theta_0, \right),
     \label{orb2}
\end{equation}
where $p$ and $q$ are positive integers and $\theta_0$ is a constant. We see that the iteration of this initial orbit has simple dynamics: $x_i^j \rightarrow x_j^{i+1}=x_{j+1}^{i}$, $y_i^j =y_j^{i+1}$, $\kappa^j=0$ and $\theta^j=\theta_0$, where the variables $x$ and $\theta$ are defined in the interval $[0,2\pi)$. The integrable periodic orbit (\ref{orb2}) is the initial point of our asymptotic procedure. The small parameter in our method is $\gamma$ and the variable $\kappa^j$ has to be small as well. We also set  $\kappa^j$ to be of order $\mathcal{O}(\gamma)$.

%|||||||||||||||||||||||||||||||||||||||||||||||||||||
\section{\large{Normal forms and sequential periodic orbits.}}
%|||||||||||||||||||||||||||||||||||||||||||||||||||||

Normal forms are reduced systems obtained by  changes of variables designed to simplify 
a set of equations \cite{Meyer92, Haragus11}. The concept of normal 
form goes back to Poincar\'e 
and for this reason, there are many different 
approaches for different problems \cite{Delshams}. 
{ When the system is written at different scales of a small parameter, it is possible to find relations between the systems dynamics and the parameters.
 These relations can also be found by performing a high number of simulations and performing parameters correlations, but this can be very time consuming and often not as accurate as the results of an asymptotic method.}

The normal form approach that we are interested { in,} is the one used in \cite{Olvera01}. 
The proposed change of variables is written as an infinite series such that the initial 
evolution equations are rewritten as a Poincar\'e-Lindstedt system ordered by a small parameter. 
In our case, the small parameter corresponds to the initial value of $\epsilon=\kappa^0$, 
where $\gamma$ is also assumed small. 
Of course, this choice is arbitrary since we are looking for periodic orbits in the map.
Therefore we can always take a different choice, for instance $\kappa^1$.

A big difference with \cite{Olvera01} is that instead of introducing a single change 
of variables for a $2$-dimensional map, for the $2N+2$-dimensional map, a change of 
variables is proposed for each { oscillator} $(x_k,y_k)$ with the \emph{a priori} 
condition that all the { oscillators} and the mean field  $(\kappa,\theta)$ have the same 
\emph{rational} winding number $\omega=2\pi p/q$, $p/q\,\in \mathbb{Q}$. The reason of 
this  condition is that we decrease the effective number of variables needed and we avoid 
performing an extra change of variables for the \emph{mean field} map.
In order to prepare the map (\ref{mapaut}) in a suitable form, we rewrite the { oscillators'}
equations (\ref{mapaut1}) and (\ref{mapaut2}) in the Lagrangian form,
\begin{subequations}
\label{maplag}
\begin{align}
  x_{k}^{n+1} -2 x_{k}^{n} +x_{k}^{n-1} &= - \kappa^{n+1}\, \sin(x_{k}^{n}-\theta^{n}) \label{maplag12}\\
  \kappa^{n+1}&= \sqrt{( \kappa^{n})^{2}+ (\eta^{n})^2} + \eta^{n} \label{maplag3}\\
  \theta^{n+1}&= \theta^{n} - \Omega + {\frac{1}{\kappa^{n+1}} \frac{\partial \eta^n}{\partial \theta^{n}}}
 \label{maplag4}
\end{align}
\end{subequations}

 The proposed change of variables for each { oscillator} is
\begin{equation}
  x_k=\zeta + \omega k + g(\zeta + \omega k,\kappa^0), \label{cambcoor}
\end{equation}
where $\omega$ is a given winding number and $g(\phi,\varepsilon)$ is given by
 \begin{equation}
    g(\phi,\varepsilon) = \sum_{j=0}^{\infty} \varepsilon^j \sum_{m\in \mathbb{Z}} g_{j,m} \,e^{ i m\phi}. \label{g_series}
 \end{equation}

Substituting (\ref{cambcoor})-(\ref{g_series}) in (\ref{maplag12}) yields for each { oscillator} the \emph{same} homological equation\footnote{The equation still depends on the label $k$ of each { oscillator}, but it can be absorbed in the angular variable: $\zeta_k:=\zeta+\omega k$}: 
\begin{equation}
 \sum_{j=0}^{\infty} (\kappa^n)^j \sum_{m\in \mathbb{Z}} g_{j,m} \,e^{ i m\zeta_k^n} c_m = \frac{\sqrt{(\kappa^n)^2+(\eta^n)^2}+\eta^n}{2i} \left\lbrace \exp\big(i (x_k^n-\theta^n)\big) - \exp\big(-i (x_k ^n-\theta^n)\big)\right\rbrace, \label{homol_gen}
\end{equation}
where 
\begin{equation}
c_m:=2(1-\cos(2\pi mp/q)), \label{c_def}
\end{equation}
\begin{equation}
\eta^n = \frac{\gamma}{2}\sum_{k'=1}^{N} \left\lbrace \exp\big(i (x_{k'}^n-\theta^n)\big) + \exp\big(-i (x_{k'}^n-\theta^n)\big)\right\rbrace. \label{eta_g}
\end{equation}

To solve (\ref{homol_gen}), we impose a relation between the two \emph{small} parameters
\footnote{For the perturbation analysis to work we can assume that both parameters 
are small. Notice that this assumption is consistent since both parameters can 
be traced back to the perturbation parameter of the standard map.} $\gamma$ and $\kappa^0$, 
so the homological equation can be organized hierarchically in terms of a single variable 
$\varepsilon=\kappa^0$, as in the Poincaré-Lindstedt method. This relation is added in order 
to have a well ordered set of equations for the perturbation analysis, and it is not a 
physical constraint that the system must necessarily obey. The lowest order at which the problem 
can be solved is: $\gamma\propto \kappa^0$ or $\gamma =\alpha \kappa^0$, $\alpha\in\mathbb{R}$. 
Also, for (\ref{eta_g}) to be summable at order $\mathcal{O}(1)$, we set the mean value 
to zero on each change of variables.

From (\ref{c_def}), it is obvious that there will be cases were (\ref{homol_gen}) will not be solvable. The terms on the right hand side that can not be eliminated, give the \emph{normal resonant form} of (\ref{maplag12})\footnote{Were both equations have been merged into 1, the \emph{lagrangian} representation of a $2$-dimensional map.}. In general the normal form has the following structure,
\be
  x_k^{n+1}-2x_k^{n}+x_k^{n-1} = \sum_{\ell=q}^{\infty} \left(\kappa^0\right)^\ell\,\left\lbrace c_\ell^{+}(\alpha)\,e^{i\ell \left(\zeta_k^n-\theta^n\right)}+ c_\ell^{-}(\alpha)\,e^{-i\ell \left(\zeta_k^n-\theta^n\right)} \right\rbrace.
\ee

The next step is to substitute the computed $g(\phi, \varepsilon)$ in the map to find the needed initial conditions and additional parameters  to have these $p/q$-\emph{periodic orbits}, that we will call \emph{sequential  periodic orbits}. 

After such substitution on the map and its next $q$ iterations, we found that for given $\kappa^0$ and $\alpha$ small, the \emph{sequential  periodic orbits} exist only for a certain value of the parameter $\Omega=\Omega(\kappa^0,\alpha)$, if the initial conditions of the { oscillators} are taken \emph{near} the \emph{fixed points}\footnote{Actually fixed points for the map iterated $q$ times.}: hyperbolic, elliptic or a mixed type.

To conclude this section we present the steps to compute the first terms of the change of variables $g(\phi,\varepsilon)$ and the resonant normal form for the sequential periodic orbit with period $\tau=2\pi /3$. We present the calculations done with them to determine the corresponding  dependence of $\Omega$ on $\kappa^0$.

%###########################################
\begin{enumerate}
  %<<<<<<<<<<<<<<<<<<<
  \item We write the zeroth order  of the homological equation (\ref{homol_gen}), substituting $(\gamma=\alpha\,\kappa^0)$\footnote{It is important to remind that the upper index $0$ is just a reference, the value of $\gamma$ will not change from iteration to iteration even though the change of variables $g$ and the \emph{normal form} are used to evaluate iterations of the map.},
\begin{equation}
 \sum_{m\in \mathbb{Z}} g_{0,m} \,e^{ i m\zeta^0} c_m = 0, \label{homol_0}
\end{equation}
where for this case is,
\begin{equation}
 c_k= \left\lbrace\begin{array}{l}
  0,\;\;\text{ k=3m} \\ 3,\;\;\text{ k=3m+1} \\ 3,\;\;\text{ k=3m+2}
 \end{array} \right. , \;\;m\in\mathbb{Z}
\end{equation}
and we have omitted the sub-index $k$ on $\zeta_k^n$ since the value of $\zeta$ does not change with $k$ in this case.
We have several free parameters that can be used later, but for the moment we set them all to zero, meaning
\begin{equation}
   g_{0,m}=0 \qquad \forall m\in \mathbb{Z}.
\end{equation}
  %<<<<<<<<<<<<<<<<<<<
 \item Now we write the next order, $(\kappa^0)$ of the homological equation,
\begin{equation}
 \sum_{k\in \mathbb{Z}} g_{1,k} e^{ i k \zeta^0}  c_k= \frac{1}{2i} \left( e^{ i (\zeta^0 - \theta^0)} - e^{-  i (\zeta^0 - \theta^0)} \right).
 \label{G1}
\end{equation}
From which we obtain,
\begin{equation}
  g_{1,1}=\frac{1}{6i}e^{-i\theta^0},\qquad 
  g_{1,-1}=\frac{-1}{6i}e^{i\theta^0}.
  \label{g13}
\end{equation}
and set the remaining $g_{1,m}$ to 0.

 %<<<<<<<<<<<<<<<<<<<
 \item We repeat the process for the following orders. 
 
       Up to order $\left((\kappa^0)^3\right)$, we compute,
 \begin{equation}
  g_{2,2}=\frac{1}{36i}e^{-i2\theta^0},\;\;\;\;\;\;\;g_{2,-2}=\frac{-1}{36i}e^{i2\theta^0},
  \label{g2n}
\end{equation}
 \begin{equation}
 g_{3,\pm 4}= -\frac{\alpha}{{96}} e^{\mp 4i\theta^0}, \;\;\;\; g_{3,\pm 2}= \frac{\alpha}{{96}} e^{\mp 2i\theta^0},\;\;\;\;
 g_{3,\pm 1}= \mp \frac{1}{432 i} e^{\mp i\theta^0}.    \label{g3n}
\end{equation}
We stop at this order due to the existence of the resonant terms,
\begin{equation}
 c_3\Big|_{0} \;g_{3,\pm 3}= \pm \frac{1}{48i} e^{\mp 3i\theta^0}
  \label{reso3}
\end{equation}

 %<<<<<<<<<<<<<<<<<<<
 \item Adding these results gives the change of variables\footnote{Where: $\phi^k:= \zeta^k - \theta^k$.}
  \begin{eqnarray}
       g(\zeta^0,\kappa^0)&=& \frac{\kappa^0}{3} \sin(\phi^0) +
                              \frac{(\kappa^0)^2}{18} \sin(2\phi^0) +
                         \frac{(\kappa^0)^3}{216} \sin(\phi^0) \nonumber \\
        && +\frac{\alpha (\kappa^0)^3}{48} \big( {\rm{cos}}(2\phi^0) -{\rm{cos}}
          (4\phi^0) \big)+\mathcal{O}((\kappa^0)^4), \label{g1_3}
  \end{eqnarray}
 
       and the resonant normal form,
  \begin{equation}
       \zeta^{1}-2\zeta^0 +\zeta^{-1} = -\frac{(\kappa^0)^3}{24} \sin(3\phi^0) 
                                        + \mathcal{O}((\kappa^0)^4).
     \label{fnr}
  \end{equation}

 %<<<<<<<<<<<<<<<<<<<
 \item Substituting (\ref{g1_3}) and (\ref{fnr}) in the map (\ref{maplag}), we solve the mean field, in other words, we look for the needed initial conditions and the parameters to have the same period for the variables $\kappa$ and $\theta$.
 \begin{enumerate}
   %>>>>>>>Kappa>>>>>>>>>
   \item ${\kappa^n}$: Substituting (\ref{g1_3}) into map in (\ref{maplag3}), we obtain,
    \begin{equation}
          \kappa^{1}= \kappa^0 \Big( 1+ \frac{\alpha (\kappa^0)^2}{6} {\sin}(3\phi^0) + O((\kappa^0)^4) \Big),
     \label{kappa1_13}
     \end{equation}
     and substituting in the following iterates yields,
     \begin{equation}
        \kappa^{3} = \kappa^0 + \frac{\alpha (\kappa^0)^3}{6} 
                     \Big( \sin(3\phi^0) + \sin(3\phi^1) + \sin(3\phi^2) \Big)
                    + \mathcal{O}((\kappa^0)^5)  \label{kappa3_13}
     \end{equation}   
     So, the condition  to have period $3$ in $\kappa^n$ is,
     \begin{equation}
       \sin(3\phi^0) + \sin(3\phi^1) + \sin(3\phi^2)
                    = \mathcal{O}((\kappa^0)^2).  \label{kappa_c}
     \end{equation}   
         
    %>>>>>>>X,Y>>>>>>>>>
   \item ${(x^n,y^n)}$: After substituting (\ref{fnr}) into (\ref{maplag12}) and  imposing $\zeta_{3}-\zeta_{0} -2\pi=0$, we obtain a similar condition,
        \begin{equation}
             \sin(3\phi^0) + \sin(3\phi^1) + \sin(3\phi^2)
                    = \mathcal{O}(\kappa^0). \label{zeta_c}
        \end{equation} 
         
   %>>>>>>>\theta>>>>>>>>>
   \item ${\theta^n}$: Substituting (\ref{g1_3}) and the results for $\kappa^n$ into the map in (\ref{mapaut4}), yields a different condition\footnote{Where ${\rm s}(\phi)\equiv\sin(\phi)$ and ${\rm c}(\phi)\equiv\cos(\phi)$.},
     \begin{eqnarray}
           \theta^3 &=& \theta^0 -3\Omega + \frac{3}{2}\alpha\kappa^0 - \frac{\alpha(\kappa^0)^2}{8}\Big[ \sum_{j=0}^{2} {\rm{c}}(3\phi_j) +
                         \frac{1}{2}\sum_{j=0}^{2} {\rm{s}}(3\phi_j)\Big] \nonumber \\
                    &&   + \frac{3\alpha(\kappa^0)^3}{324} + \frac{\alpha^2(\kappa^0)^3}{24} [ {4 \rm{s}}(3\phi^0) + 3 {\rm{s}}(3\phi^1) 
                         - {\rm{s}}(3\phi^2) ]+ \mathcal{O}(\alpha(\kappa^0)^4). \nonumber \\   
         \label{theta_c}
         \end{eqnarray}   
   
  \end{enumerate} 
  
   The conditions (\ref{kappa_c}), (\ref{zeta_c}) and (\ref{theta_c}) can only be solved in the general case  if we restrict to the fixed points {(of the map iterated $q$ times)}. After the change of variables, $\phi=\frac{2n\pi}{3}$ ($\phi=\frac{ (2n+1)\pi}{3}$) corresponds to the linear elliptic (hyperbolic) {fixed} points at least at order $\kappa^0$.
   
   Then (\ref{kappa_c}) and (\ref{zeta_c}) are satisfied and (\ref{theta_c}) yields,
       \begin{equation}
          \Omega = \frac{\alpha\kappa^0}{2} - \frac{\alpha(\kappa^0)^2}{8} +
                   \frac{\alpha(\kappa^0)^3}{324} + \mathcal{O}(\alpha (\kappa^0)^4),
          \label{omelip}
       \end{equation}
  for linear elliptic \emph{fixed} points and 
       \begin{equation}
          \Omega = \frac{\alpha\kappa^0}{2} + \frac{\alpha(\kappa^0)^2}{8} +
                   \frac{\alpha(\kappa^0)^3}{324} + \mathcal{O}(\alpha (\kappa^0)^4),
          \label{omhip}
       \end{equation}
     for linear hyperbolic \emph{fixed} points.
     
     If $\kappa^1=\sqrt{(\kappa^0)^2+(\eta^0)^2}+\eta^0\simeq \kappa^0+\eta^0+\mathcal{O}((\eta^0)^2/\kappa^0)$, we can expand $\eta^0=\gamma\sum\sin(x_k^0-\theta^0)$ and we substitute $x_k^0$ from (\ref{cambcoor}), we obtain
     \begin{equation}
          \kappa^1-\kappa^0\simeq \frac{\alpha^2 (\kappa^0)^4}{64}.
          \label{dkappa}
     \end{equation}
\end{enumerate}
%###########################################
\vspace{-0.3cm}
%////////////////////////////
  \begin{figure}[h!]
     \centering
    \includegraphics[height=9.3cm]{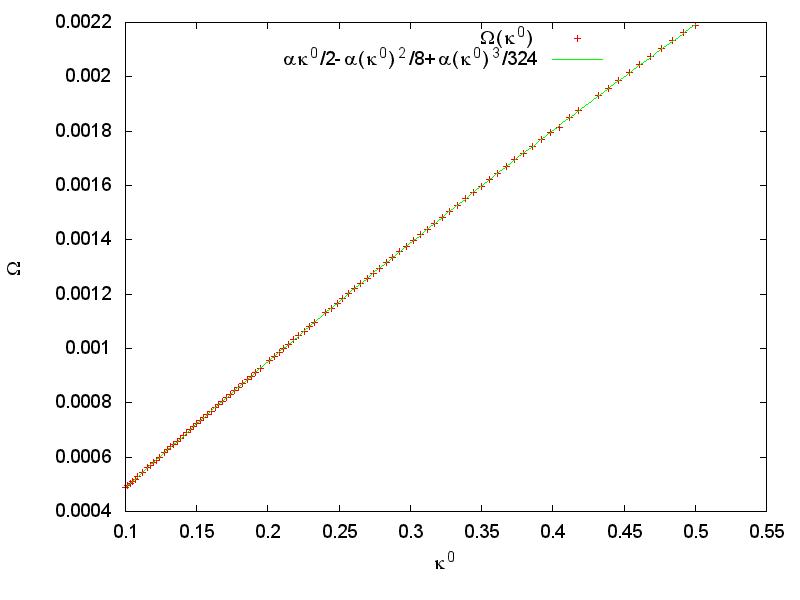}
     \caption{\footnotesize{The figure displays {in red} the values of $\Omega$ for a given $\kappa^0$ (with $\alpha=\frac{1}{100}$) found by using numerical continuation of a  sequential {periodic} orbit with winding number $\omega=\frac{2\pi}{3}$, that started from the linear elliptic \emph{fixed} points. { The overlapped green line correspond to parameter relation (\ref{omelip}) found with normal forms.}}}
     \label{fig_elip}
 \end{figure}
%////////////////////////////
 
%////////////////////////////
  \begin{figure}[h!]
     \centering
     \includegraphics[height=9.3cm]{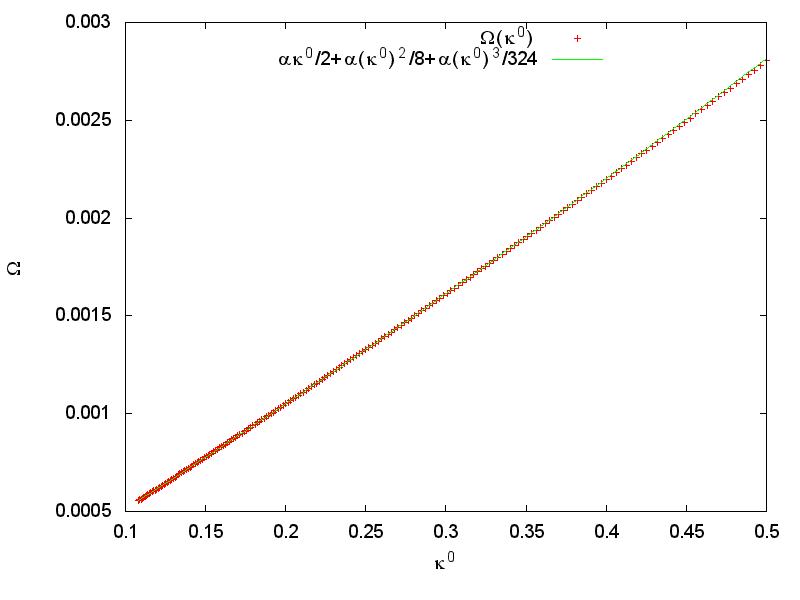}
     \caption{\footnotesize{The figure displays  in red the values of $\Omega$  for a given $\kappa^0$ (with $\alpha=\frac{1}{100}$) found by using numerical continuation of a   sequential periodic orbit with winding number $\omega=\frac{2\pi}{3}$, that started from the linear hyperbolic \emph{fixed} points. { The overlapped green line correspond to parameter relation (\ref{omhip}) found with normal forms.}}}
     \label{fig_hyp}
 \end{figure}
%////////////////////////////

Table \ref{tabla1} shows the change of variables, normal resonant form and $\Omega(\kappa^0)$ relations for two kinds of initial conditions (elliptic and hyperbolic fixed points) for given rotation numbers.

\begin{table}[h!]
\centering
\begin{tabular}{|c|l|}\hline
$1/q$  &  \\  \hline \hline
 \multirow{10}{0.6cm}{\centering 1/2} &\\
    &  $ g(\zeta^0,\kappa^0)= \frac{\kappa^0}{4} \sin(\phi^0) +
        \frac{\alpha(\kappa^0)^2}{32} \left( \cos(\phi^0)- \cos(3\phi^0)\right)
        - \frac{(3-\alpha^2)(\kappa^0)^3}{512} \sin(\phi^0)$\\
    &  \phantom{OOOOOO}  $ + \frac{(1+\alpha^2)(\kappa^0)^3}{512} \sin(3\phi^0) 
        - \frac{\alpha^2(\kappa^0)^3}{512} \sin(5\phi^0) +  \mathcal{O}((\kappa^0)^4)$\\
    & \\
    & $\zeta^{1}-2\zeta^0 +\zeta^{-1} = -\frac{(\kappa^0)^2}{a} \sin(2\phi^0) 
               + \frac{\alpha(\kappa^0)^3}{32}  \left( \cos(4\phi^0)-1 \right) 
             +\mathcal{O}((\kappa^0)^4)$  \\ 
    & \\
        &  \emph{Elliptic:}\phantom{Oii} $\Omega = \mathcal{O}((\kappa^0)^4)$\\ &\\
    &  \emph{Hyperbolic:} $\Omega = \frac{\alpha\kappa^0}{2} -
           \frac{\alpha(\kappa^0)^2}{64} (1+\alpha^2)+ \frac{3\alpha(\kappa^0)^3}{128} 
           + \mathcal{O}(\alpha (\kappa^0)^4) $\\ &\\ \hline
%------------------------------------
 \multirow{10}{0.6cm}{\centering 1/3} &\\
    &  $ g(\zeta^0,\kappa^0)= \frac{\kappa^0}{3} \sin(\phi^0) +
                              \frac{(\kappa^0)^2}{18} \sin(2\phi^0) 
                               + \frac{(\kappa^0)^3}{216} \sin(\phi^0)
                                +\frac{\alpha (\kappa^0)^3}{48} \big( \cos(2\phi^0) 
                     -\cos  (4\phi^0) \big)$  \\
    & \phantom{OOOOOO} $     +\mathcal{O}((\kappa^0)^4)$  \\ &\\

    &  $\zeta^{1}-2\zeta^0 +\zeta^{-1} = -\frac{(\kappa^0)^3}{24} \sin(3\phi^0) 
                 + \mathcal{O}((\kappa^0)^4)$  \\ 
    & \\
    &  \emph{Elliptic:}\phantom{Oii} $\Omega = \frac{\alpha\kappa^0}{2} -
                \frac{\alpha(\kappa^0)^2}{8} + \frac{\alpha(\kappa^0)^3}{324} +
                 \mathcal{O}(\alpha (\kappa^0)^4) $\\  &\\
    &  \emph{Hyperbolic:} $\Omega = \frac{\alpha\kappa^0}{2} -
           \frac{\alpha(\kappa^0)^2}{8} + \frac{\alpha(\kappa^0)^3}{324} 
           + \mathcal{O}(\alpha (\kappa^0)^4) $\\ 
    &\\ \hline
%------------------------------------
 \multirow{13}{0.6cm}{\centering 1/6} &\\
    &  $g(\zeta^0,\kappa^0)= \kappa^0 \sin(\phi^0) 
             + \frac{(\kappa^0)^2}{6} \sin(2\phi^0) 
             - \frac{7(\kappa^0)^3}{24} \sin(3\phi^0)
             + \frac{(\kappa^0)^3}{96} \sin(3\phi^0)
             - \frac{43(\kappa^0)^4}{576} \sin(2\phi^0)$ \\
    & \phantom{OOOOOO}$ + \frac{17(\kappa^0)^4}{576} \sin(4\phi^0)
             + \frac{25(\kappa^0)^5}{288} \sin(\phi^0) 
             - \frac{185(\kappa^0)^5}{1152} \sin(3\phi^0)
             + \frac{51(\kappa^0)^5}{384} \sin(5\phi^0)$ \\ 
    & \phantom{OOOOOO}$ + \frac{211(\kappa^0)^6}{864} \sin(2\phi^0) 
             - \frac{4227(\kappa^0)^6}{34560} \sin(4\phi^0)
             + \frac{1077\alpha(\kappa^0)^6}{3840} 
                 \left( \cos(7\phi^0)-\cos(5\phi^0) \right)$ \\
    & \phantom{OOOOOO}$  + \mathcal{O}((\kappa^0)^7) $ \\&\\
    &  $\zeta^{1}-2\zeta^0 +\zeta^{-1} = -\frac{1077(\kappa^0)^6}{3840} \sin(6\phi^0) 
                       + \mathcal{O}((\kappa^0)^7)$  \\ 
    & \\
    &  \emph{Elliptic:}\phantom{Oii} $\Omega = \alpha\kappa^0 
                +\frac{9}{24}\alpha (\kappa^0)^3 
                -\frac{379}{1151} \alpha (\kappa^0)^5 
                +\frac{1077}{3840} \alpha (\kappa^0)^5 +\mathcal{O}((\kappa^0)^7)$\\ &\\
    &  \emph{Hyperbolic:} $\Omega = \alpha\kappa^0 
                +\frac{9}{24}\alpha (\kappa^0)^3 
                -\frac{379}{1151} \alpha (\kappa^0)^5 
                -\frac{1077}{3840} \alpha (\kappa^0)^5 +\mathcal{O}((\kappa^0)^7)$\\ &\\ \hline
\end{tabular}
\caption{\footnotesize{Normal forms for some sequential periodic orbits with ration rotation numbers. { In the table we display: (i) the change of variables $g$, (ii) the resonant normal form and needed parameter relations for (iii)\emph{elliptic} and (iv)\emph{hyperbolic} fixed points.}}}
\label{tabla1}
\end{table}

The three examples show that the relation of the asymptotic value of the mean field
variables with the parameters of the self--consistent map (\ref{mapaut}). For each
rotation number of the sequential periodic orbit we obtain a specific value of the
average of $\kappa$ and the amplitude of its oscillation. Table \ref{tabla1}
shows these values for rotation numbers $1/2$, $1/3$ and $1/6$. The main point
is the relation between $\Omega$ and $\kappa^0$.  Fixing  the values of the
parameters $\gamma$ and $\Omega$, we can find the average value of kappa, where
$\bar{\kappa} = \kappa^0$ by using equations (\ref{omelip}) and (\ref{omhip}) and
the relation $\gamma=\alpha\kappa^0$. The amplitude of the variation of
$\kappa$ is  given by $\Delta\,\kappa= {\rm max} | \kappa^i - \kappa^0 |$, for
$i=1,\dots,N$, which first approximation is written in (\ref{dkappa}). This method can be used for any sequential periodic orbits of period
$\tau$.

 We use our numerical procedure to compute the previous sequential  periodic orbits.
The first step is to compute a periodic orbit with the same rotation number as
the sequential periodic orbit of the standard map (\ref{stnmap}) with value of
the parameter $\kappa= \kappa^0$. This orbit was our initial guess for the
numerical approximation of the sequential periodic orbit. In order to have
convergence in the Newton method,  it is  necessary to define $\Omega$ as a free
parameter, in this way the numerical procedure converges for a specific
value of $\Omega$.

The normal forms and the numerical method to find sequential  periodic orbits imply an
interesting relation between $\kappa^0$ and $\Omega$.  In order to compare the
numerical and the asymptotical methods, we evaluate the relation (\ref{omelip})
using the values of $\kappa$ and $\Omega$ obtained from the numerical method.
Figure \ref{fig_elip} (resp. \ref{fig_hyp}) shows the evaluation of (\ref{omelip}) (resp. (\ref{omhip})) using the numerical
values of $\kappa$ versus the numerical value of $\Omega$, the figures show a
very good match of the numerical procedure and the normal forms. In this form
we obtain very good agreement of the numerics and the asymptotic procedures
using normal forms for the case of sequential periodic orbits.

The numerical computations shows that $\Delta\kappa\simeq\frac{\alpha^2}{40}(\kappa^0)^4$ for linear hyperbolic fixed points, which is of the order of the asymptotic result $\Delta\kappa\simeq\frac{\alpha^2}{64}(\kappa^0)^4$ from (\ref{dkappa}).

{ A conclusion that we obtain of the use of normal forms is that the sequential periodic orbits that we have constructed, have associated values $\bar{\kappa}$ and $\Delta\kappa$ as functions of the parameters $\Omega$ and $\gamma$. This is a possible explanation of why the mean field variable $\kappa$ achieves a mean value and oscillates, at least for this sequential periodic orbit scenario.}

{ In a more general dynamic, we could picture that the oscillatory evolution of the variable $\kappa$ is driven by a set of orbits close to periodic orbits, periodic orbits associated to the parameters of the evolution. Due this hypothesis, we are interested in studying the frequency space of $\kappa$ for long times. We hope to perform this analysis in a future work.}

%|||||||||||||||||||||||||||||||||||||||||||||||||||||
\section{\large{Conclusions}}
%|||||||||||||||||||||||||||||||||||||||||||||||||||||

In this paper we studied self-consistent chaotic transport in a mean-field Hamiltonian model.
The model consists of $N \gg 1$ coupled area-preserving twist maps in which the amplitude and phase of the perturbation (rather than being constant like in the standard map) are dynamical variables. 
The model provides a simplified description of transport in marginally stable systems including
vorticity mixing in strong shear flows and electron dynamics in plasmas.  

Of particular interest was the study of  coherent structures and periodic orbits.
Numerical simulations showed that self-consistency leads to the formation of a coherent macro-particle trapped around 
the elliptic fixed point of the system, accompanied by an oscillatory  behavior of the mean field. To study in detail 
the transport properties of the self-consistent map in this case, we introduced a 
non-autonomous standard map in which the amplitude of the perturbation alternates between two values, $\kappa_1$ and 
$\kappa_2$ mimicking the observed oscillatory behavior of the mean field in the full self-consistent map.  

A  turnstile-type  mechanism that allows transport across instantaneous KAM invariant circles in non-autonomous systems was presented. This mechanism explains how, contrary to intuition, orbits in the 
non-autonomous standard map can cross the golden-mean torus even though the values of $\kappa_1$ and $\kappa_2$ might be below  the critical threshold for the destruction of the golden mean torus in the standard map $\kappa_c$. This result  raises the question of the critical  $\kappa_1$ and $\kappa_2$ parameter values for the onset of global transport (i.e. complete absence of KAM invariant circles) in the non-autonomous standard map. We approached this question numerically and found the region in the
$(\kappa_1,\kappa_2)$ for which global transport exits. Interestingly, it was observed that the boundary of this ``horn-shaped" region exhibits a cusp singularity at $\kappa_1=\kappa_2=\kappa_*$. 

Based on Moser's twist theorem,  a study of the twist condition in the non-autonomous standard  map in Eq.~(\ref{map}) indicates that, for small enough values of $\kappa_1$ and $\kappa_2$
the map must have invariant circles. That is,  there will not be diffusion no matter 
how many iterations are performed on the map. Moreover, since the same argument can be 
given for a generalized version of the map (\ref{map}) with $m$ parameters $\kappa_i$, 
$i=1,2...,m$, it is possible to construct a high dimensional system with no diffusion.

We also showed that a type of periodic
orbits, refer to as sequential periodic orbits,  can be used to explain the
behavior of the mean field variables of the self-consistent map. 
In particular, using normal forms,  we showed that for sequential periodic orbits to exists, a specific 
relationship between the mean-field variables
$\kappa^0$ and $\theta^0$  and the map parameters
$\gamma_k$ and $\Omega$ must hold. 
This result opens the possibility of predicting the asymptotic
values of  mean-field variable in the case when  the initial condition corresponds to a
sequential periodic orbit. { Nevertheless, at this moment we are incapable to establish a direct connection 
between these periodic orbits and the oscillatory behavior observed in figures \ref{catseye_map} and \ref{evol}.}
 
Among the several problems we plan to explore in the future is to determine if one can find periodic orbits that drive the asymptotic dynamics of the mean field variables.
In this way, we would be able to estimate 
the transport of the self-consistent model with the help of the
turnstile-type mechanism.
Another interesting problem is the study the self-consistent map 
in the limit of small $N$. Preliminary results indicate that, using symmetries and conservation laws, for $N=1$,  the $4$-dimensional map in Eq.~(\ref{mapaut})
can be reduced to a $2$-dimensional map.

%|||||||||||||||||||||||||||||||||||||||||||||||||||||
\section*{\large{Acknowledgments}}
%|||||||||||||||||||||||||||||||||||||||||||||||||||||
This work was founded by PAPIIT IN104514, FENOMEC-UNAM and
by the Office of Fusion Energy Sciences of the US Department
of Energy at Oak Ridge National Laboratory, managed by UT-Battelle, LLC, for the
U.S.Department of Energy under contract DE-AC05-00OR22725. { We also express our gratitude  to the graduate program in Mathematics of UNAM for making the GPU servers available to perform our computations and especially to Ana Perez for her invaluable help. Finally, we would also like to thank the anonymous referee whose valuable comments have improved the 
presentation of the paper.}

\end{document}